\documentclass[a4paper, 11pt]{amsart}
\usepackage{amsfonts}
\usepackage{amsmath}
\usepackage{tikz-cd}
\usepackage{mathtools}
\usepackage{amssymb}
\usepackage{amsthm}
\usepackage{amscd}
\usepackage{color}

\usepackage[colorlinks]{hyperref}
\definecolor{linkblue}{RGB}{1,1,190}
\definecolor{citered}{RGB}{190,1,1}  
\hypersetup{
	linkcolor=linkblue,
	urlcolor=linkblue,
	citecolor=citered
}

\newcommand{\N}{\mathbb N}

\DeclareMathOperator{\Pic}{Pic}

\DeclareMathOperator{\spec}{spec}

\def\P{\ensuremath\mathfrak{P}}
\def\p{\ensuremath\mathfrak{p}}
\def\O{\ensuremath\mathcal{O}}
\def\f{\ensuremath\mathfrak{f}}

\begin{document}

\title[A characterization of half-factorial orders in number fields]{A characterization of half-factorial orders \\ in algebraic number fields} 

\author{Balint Rago}

\address{University of Graz, NAWI Graz, Department of Mathematics and Scientific Computing, Heinrichstraße 36,
8010 Graz, Austria}

\thanks{This work was supported by the Austrian Science Fund FWF, Project Number W1230}

\email{balint.rago@uni-graz.at}



\begin{abstract}
    We give an algebraic characterization of half-factorial orders in algebraic number fields. This generalizes prior results for seminormal orders and for orders in quadratic number fields.
    \end{abstract}

\maketitle

\smallskip
\section{Introduction}
\smallskip

\noindent Every nonzero nonunit of a Noetherian domain $R$ can be written as a finite product of irreducible elements. The domain $R$ is said to be half-factorial if, for every such $a \in R$ and for each two factorizations into irreducibles $a = u_1 \cdots u_k = v_1 \cdots v_{\ell}$, it follows that $k = \ell$. Thus, every factorial domain is half-factorial.

In $1960$, Carlitz \cite{Ca60} proved that the principal order of an algebraic number field is half-factorial if and only if its class group has at most two elements.
Ever since then,
half-factoriality has been a central topic in the factorization theory of monoids and domains (we can refer only to some of the many contributions, see \cite{Ch-Co00, Co03, Co05a, Sc05c,Pl-Sc05a, Pl-Sc05b, Ma-Ok09a, Ro11a,  Ma-Ok16a, G-L-T-Z21}). In particular, orders in algebraic number fields received constant attention.  In $1983$,  Halter-Koch gave a  characterization of half-factorial orders in quadratic number fields (\cite{HK83a}, \cite[Theorem 3.7.15]{Ge-HK06a}), and a characterization of half-factoriality for seminormal orders was given in \cite[Theorem 6.2]{Ge-Ka-Re15a}. In spite of further partial results \cite{Pi00, Pi02a, Ka05b, Ph12b}, the case of arbitrary orders in algebraic number fields remained open.

The main result of the present paper establishes an algebraic characterization of half-factoriality for arbitrary orders in algebraic number fields. The known characterizations of half-factorial seminormal orders and of half-factorial orders in quadratic number fields are obtained as corollaries of our main result (Corollaries \hyperref[4.6]{4.6.} and \hyperref[4.7]{4.7.}). \\

\smallskip

\phantomsection
\label{1.1}\noindent \textbf{Theorem 1.1.} \textit{Let $K$ be an algebraic number field with principal order $\O_K$, and let $\O \subseteq \mathcal O_K$ be an order in $K$ with conductor $\f=\P_1^{k_1}\ldots\P_s^{k_s}$, where $s, k_1, \ldots, k_s \in \N$, $\P_1, \ldots, \P_s \in \spec (\mathcal O_K)$, and we set  $\p_i=\P_i\cap\O$. Then $\O$ is half-factorial if and only if the following conditions are satisfied.} \newpage

\begin{itemize}

\item[(i)] $|\Pic (\O_K)| \le 2$. \\ 

\item[(ii)] $\O\cdot\O_K^\times=\O_K$. \\

\item[(iii)] \textit{For all $i \in [1,s]$, we have $k_i\leq 4$ and $\mathsf{v}_{p_i}(\mathcal{A}(\O_{\p_i}))\subseteq\{1,2\},$ where $p_i$ is an arbitrary prime element of $(\O_K)_{\p_i}$. If $\P_i$ is principal, then $k_i\leq 2$ and $\mathsf{v}_{p_i}(\mathcal{A}(\O_{\p_i}))=\{1\}$.} \\
\end{itemize}

It is well-known that condition (i) is equivalent to $\O_K$ being half-factorial and that condition (iii) for principal ideals is equivalent to the local order $\O_{\p_i}$ being half-factorial (see Proposition \hyperref[3.2]{3.2.}). At the end of section 4, we will present an example of a half-factorial order $\O$ with the property $\mathsf{v}_{p_i}(\mathcal{A}(\O_{\p_i}))=\{1,2\}$ for some prime ideal $\p_i$, in order to show that none of the conditions in the theorem is vacuous. \\

In Section 2 we give a quick overview of some basic properties of orders in algebraic number fields and introduce the notion of irreducible orders, which allows the decomposition of an arbitrary order into simpler components. Furthermore, we give a short introduction to factorization theory in multiplicative monoids.  \\

In Section 3 we will discuss algebraic and arithmetic properties of the local orders $\O_\p$. In particular, we will recall their characterization of half-factoriality, given in \cite{Ka05b}. \\

In Section 4 we prove our main theorem and give two corollaries, proving the aforementioned characterization of half-factoriality in seminormal (resp. quadratic) orders. \\

\smallskip
\section{Preliminaries}
\smallskip

Let $K$ be an algebraic number field with principal order $\O_K$. An order in $K$ is a subring $\O\subseteq \O_K$, with quotient field $\mathsf{q}(\O)=K$. Every order is Noetherian and one-dimensional. The ideal \[\f=\{x\in\O_K\ | \ x\O_K\subseteq\O\}\] is the greatest ideal of $\O_K$ which is contained in $\O$. It is called the conductor of $\O$. \\

Apart from the case of quadratic number fields, there is no explicit description of orders in $\mathcal O_K$. However, there is a characterization of which ideals of $\mathcal O_K$ occur as conductor ideals (\cite[Theorem 2.12.9]{HK20a}). For recent work on orders in cubic number fields and for work on the distribution of orders, we refer to \cite{Le-Pr16a, Ka-Ma-Ta15}. \\

\noindent The residue class ring $\O/\f$ is a subring of $\O_K/\f$ and we have an injection \[\{\text{Orders with conductor } \f\}\longrightarrow\{\text{Subrings of }\O_K/\f\}.\] Let \[\f=\P_1^{k_1}\cdots \P_s^{k_s},\] be the unique factorization of $\f$, where $s\in \N$ and $\P_1,\ldots,\P_s$ are pairwise distinct prime ideals of $\O_K$.
Then, by the Chinese Remainder Theorem, we have an isomorphism \[\O_K/\f\cong \O_K/\P_1^{k_1}\times\cdots\times \O_K/\P_s^{k_s}.\] Let $\p_1,\cdots,\p_r$ be the prime ideals of $\O$ that contain $\f$ and let, for $i\in [1,r]$, $\P_{i,1},\cdots,\P_{i,m_i}$ be the prime ideals of $\O_K$ that lie over $\p_i$. Then we can write \[\f=\f_1\cdots\f_r,\] where \[\f_i=\P_{i,1}^{k_{i,1}}\cdots\P_{i,m_i}^{k_{i,m_i}}\] and \[\sum_{i=1}^{r}m_i=s.\] We call $\O$ \textit{irreducible}, if $r=1$, meaning that there is only one prime ideal $\p$ containing the conductor. \\

Let $\pi_i:\O/\f\longrightarrow \O_K/\f_i$ be the restriction of the canonical projection $\O_K/\f\rightarrow \O_K/\f_i$. Since $\O$ is one-dimensional, we have $\p_j\not\subseteq \p_i$ for every $j\in [1,r]\setminus\{i\}$ and thus we find an element $x_j\in\p_j\setminus\p_i$. Then, for a sufficiently large $n\in\mathbb{N}$, we see that \[x=\prod_{j\in [1,r]\setminus\{i\}}x_j^n\in\bigcap_{j\in [1,r]\setminus\{i\}}\f_j\setminus\p_i,\] i.e. $\pi_j(x+\f)=0$ for $j\neq i$ and $\pi_i(x+\f)$ is a unit in $\O_K/\f_i$. 
This implies that there is an isomorphism \[\O/\f\cong\pi_1(\O/\f)\times\cdots\times\pi_r(\O/\f)\] and we can find orders $\O_i$ with conductor $\f_i$ respectively such that $\O_i/\f_i=\pi_i(\O/\f)$ and \[\O=\bigcap_{i=1}^{r}\O_i.\] We will call the orders $\O_i$ the \textit{irreducible components} of $\O$. \\

We will call the map \[\spec(\O_K)\longrightarrow \spec(\O)\] \[\P\mapsto\P\cap\O\] the spec-map. This map is surjective, since $\O\subseteq \O_K$ is an integral extension. It is injective if and only if every prime ideal of $\O$ has precisely one prime ideal lying over it, which is equivalent to every irreducible component of $\O$ having a conductor of the form $\P^k$ for some prime ideal $\P$ of $\O_K$. The following lemma will be useful later on. \\

\phantomsection
\label{2.1} \noindent \textbf{Lemma 2.1.} \textit{Let $\O$ be an order with conductor $\f$ and suppose that the spec-map is bijective. Then for every $x\in\O_K$, there is $n\in\mathbb{N}$ such that $x^n\in\O$.} \\

\noindent \textit{Proof:} We can identify $\O/\f$ with $\O_1/\f_1\times\cdots\times \O_r/\f_r,$ where the $\O_i$'s are the irreducible components of $\O$. Let $x\in\O_K$. Since every $\f_i$ is a power of a prime ideal of $\O_K$, the $i$-th component of $x+\f$ is either nilpotent or a unit in $\O_K/\f_i$. Hence there is $n\in\mathbb{N}$ such that every component of $x^n+\f$ is either $0$ or $1$, which implies that $x^n\in\O_i$ for all $i$. Since \[\O=\bigcap_{i=1}^{r} \O_i,\] we are done. \qed \\

\noindent We now gather some terminology from factorization theory. The main reference is \cite{Ge-HK06a}. \\

By a \textit{monoid}, we mean a commutative, cancellative semigroup with unit element. Let $H$ be a monoid. We denote by $\mathcal{A}(H)$ the set of atoms (irreducible elements) of $H$, by $H^\times$ the group of units (invertible elements) of $H$ and by $H_{\text{red}}=H/H^\times=\{aH^\times \ |\  a\in H\}$ the associated reduced monoid. In general, we call a monoid $H$ \textit{reduced} if $H^\times=\{1\}$. $H$ is called \textit{atomic} if every non-unit $a\in H$ can be written as a product of irreducible elements of $H$. Let $R$ be an integral domain. We will denote the multiplicative monoid of non-zero elements of $R$ by $R^\bullet$. All arithmetic terms defined for monoids carry over to domains and we set $\mathcal{A}(R)=\mathcal{A}(R^\bullet)$ and so on. \\

Let $H$ be an atomic monoid and suppose that $a\in H\setminus H^\times$ can be factorized as \[a=a_1\cdots a_s,\] where $a_1,\cdots,a_s\in\mathcal{A}(H)$. We call the integer $s$ the \textit{length} of the factorization and \[\mathsf{L}(a)=\{s \ |\  a=a_1\cdots a_s \text{ for irreducible } a_i\}\] the \text{set of lengths} of $a$. For a unit $\varepsilon\in H^\times$, it is a convention to set $\mathsf{L}(\varepsilon)=\{0\}$. We call \[\rho(a)=\frac{\sup \mathsf{L}(a)}{\min \mathsf{L}(a)}\in\mathbb{Q}_{\geq 1}\cup \{\infty\}\] the \textit{elasticity} of $a$ and $\rho(H)=\sup\{\rho(a)\ |\ a\in H\setminus H^\times\}$ the \textit{elasticity} of $H$. $H$ is called \textit{half-factorial} if $|\mathsf{L}(a)|=1$ for every $a\in H$. Clearly, $H$ is half-factorial if and only if $\rho(H)=1$. \\

\smallskip
\section{Local Orders}
\smallskip

In this section we will discuss algebraic and arithmetic properties of the localizations $\O_\p$ of an order $\O$ with conductor $\f$. If the prime ideal $\p$ does not contain $\f$, then $\O_\p$ is a DVR and hence factorial. If $\p\supseteq\f$, then $\O_\p$ is not integrally closed. Its integral closure is the semilocal ring \[\O_K\cdot {(\O\setminus\p)}^{-1},\] which we will denote by $(\O_K)_\p$. In the following, we will assume that the spec-map is bijective, in which case $(\O_K)_\p$ is a DVR. Let $p$ be a prime element of $(\O_K)_\p$. 
We will denote the valuation on $(\O_K)_\p$ with respect to $p$ by $\mathsf{v}_p$. Then we can write \[\O_\p=\bigcup_{i\in\mathbb{N}_0}p^i\cdot U_{i,p}\cup\{0\},\] where $U_{i,p}=\{u\in(\O_K)_\p^\times \ | \ p^iu\in\O_\p \}$. Note that for a positive integer $k$, $U_{k,p}\neq\emptyset$  is equivalent to $k\in\mathsf{v}_p(\O_\p)$. There is a smallest integer $\alpha$ such that $U_{j,p}=(\O_K)_\p^\times$ for all integers $j\geq\alpha$, which we call the \textit{exponent} of $\O_\p$. Let $\P$ be the prime ideal of $\O_K$, lying over $\p$. It is easy to check that the exponent of $\O_\p$ is the highest power of $\P$ that divides $\f$. \\

Suppose for a moment that $\O$ is irreducible, i.e. $\f=\P^k$ for some $k\in\mathbb{N}$. Then we can express the reduced monoid ${(\O_\p^\bullet)}_\text{red}$ as \[\bigcup_{i\in\mathbb{N}_0}\bar{p}^i \cdot\overline{U}_{i,\bar{p}},\] where $\bar{p}=p+\f\in\O_K/\f$ for some $p\in\P\setminus\P^2$ and \[\overline{U}_{i,\bar{p}}=\{u\in (\O_K/\f)^\times/(\O/\f)^\times \ | \ \bar{p}^i u\subseteq \O/\f\}.\]  

The following lemma shows that passing from an order to an irreducible component preserves localizations. \\

\phantomsection
\label{3.1} \noindent \textbf{Lemma 3.1.} \textit{Let $\O$ be an order with conductor $\f$ and let $\p$ be a prime ideal which contains $\f$. Let $\tilde{\O}$ be the irreducible component of $\O$ which corresponds to $\p$, let $\tilde{\f}$ be the conductor of $\tilde{\O}$ and let $\tilde{\p}$ be the unique prime ideal of $\tilde{\O}$, which contains $\tilde{\f}$. There is an isomorphism of monoids} \[{(\O_{\p}^\bullet)}_{\text{red}}\cong {(\tilde{\O}_{\tilde{\p}}^\bullet)}_{\text{red}}.\]

\noindent \textit{Proof:} The inclusions $\O\subseteq\tilde{\O}$ and $\p\subseteq\tilde{\p}$ imply that $\O_\p\subseteq\tilde{\O}_{\tilde{\p}}$ and $\O_\p^\times\subseteq\tilde{\O}_{\tilde{\p}}^\times$. This means that \[\psi:{(\O_{\p}^\bullet)}_{\text{red}}\longrightarrow{(\tilde{\O}_{\tilde{\p}}^\bullet)}_{\text{red}}\]\[x\O_\p^\times\mapsto x\tilde{\O}_{\tilde{\p}}^\times\] is a well-defined homomorphism of monoids. Suppose that $\psi(x\O_\p^\times)=\psi(y\O_\p^\times)$ for elements $x,y\in \O_\p^\bullet$. We can assume that $x,y \in \O$. 
Then $x$ and $y$ are associate in $\tilde{\O}_{\tilde{\p}},$ i.e. there are $a,b\in\tilde{\O}\setminus\tilde{\p}$ such that $ax=by$. By the Chinese Remainder Theorem, we can find an element $c\in \f\cdot\tilde{\f}^{-1}$, where $\tilde{\f}^{-1}$ denotes the fractional ideal, inverse to $\tilde{\f}$, such that $c+\tilde{\f}=1\in\O_K/\tilde{\f}$. Then $ac,bc\in\O\setminus\p$ and $x$ and $y$ are associate in $\O_\p$. Hence $\psi$ is injective. \\

For surjectivity, it is sufficient to show that every $x\in\tilde{\O}$ is associate to some $y\in\O$ in $\tilde{\O}_{\tilde{\p}}$. Again, by taking an element $c\in\f\cdot\tilde{\f}^{-1}$ with $c+\tilde{\f}=1$, we have $cx\in\O$. \qed \\

In \cite{Ka05b}, Kainrath gave a characterization of half-factorial local orders, which we will now present. Suppose again that the spec-map is bijective, so that $(\O_K)_\p$ is a DVR. Let $k$ and $\bar{k}$ be the residue fields of $\O_\p$ and $(\O_K)_\p$ respectively and let $\pi:(\O_K)_\p\rightarrow\bar{k}$ be the canonical surjection. We will identify $k$ with $\pi(\O_\p)$. Furthermore, let $V_{i,p}=\pi(U_{i,p})\cup\{0\}$, where $p$ is some prime element of $(\O_K)_\p$. \\

\phantomsection
\label{3.2} \noindent \textbf{Proposition 3.2.(\cite[Theorem 3.3.]{Ka05b})} \textit{The following statements are equivalent.} \\ 

(a) \textit{The local order $\O_\p$ is half-factorial.} \\

(b) \textit{$U_{1,p}\cdot U_{1,p}=(\O_K)_\p^\times$.} \\

(c) \textit{$V_{1,p}\cdot V_{1,p}=\bar{k}.$} \\

The proposition shows that $\O_\p$ being half-factorial is equivalent to $\mathsf{v}_p(\mathcal{A}(\O_\p))=\{1\}$, since $p\cdot U_{1,p}$ is precisely the set of elements of value $1$. Since $U_{1,p}\cdot U_{1,p}\subseteq U_{2,p}$, we immediately see that the exponent of $\O_\p$ is at most two. This explains the bound on the exponents of the conductor for principal prime ideals in condition (iii) in Theorem \hyperref[1.1]{1.1.} For the other bound we need to prove the following \\

\phantomsection
\label{3.3} \noindent \textbf{Proposition 3.3.} \textit{Suppose that the spec-map is bijective and that $\mathsf{v}_p(\mathcal{A}(\O_\p))\subseteq\{1,2\}$, where $p$ is a prime element of $(\O_K)_\p$. Then the exponent of $\O_\p$ is at most four.} \\

\noindent \textit{Proof:} Clearly, $\mathsf{v}_p(\mathcal{A}(\O_\p))=\{2\}$ is impossible since $\O_\p$ contains infinitely many elements of odd value. Hence, $\O_\p$ contains an element of value $1$ and so we can assume that $p\in\O$. We have an ascending chain \[1\in V_{1,p}\subseteq V_{2,p}\subseteq\cdots \subseteq V_{i,p}\subseteq\cdots,\] which stabilizes at some point. By \cite[Lemma 3.2.]{Ka05b}, we have $U_{i,p}=(\O_K)_\p^\times$ if and only if $V_{i,p}=\bar{k}$ and so it is enough to show that $V_{4,p}=\bar{k}$. Since every element of value 4 can be written as a product of elements of value 2, we obtain $V_{4,p}\subseteq V_{2,p}\cdot V_{2,p}$ and since the reverse inclusion trivially holds, we have equality. 
Since $V_{4,p}$ is a $k$-subspace of $\bar{k}$, it is also a field by \cite[Theorem 2.6.]{Ka05b}. Similarly, we have  $V_{8,p}= V_{4,p}\cdot V_{4,p}$, which then implies that $V_{8,p}=V_{4,p}$ and by induction, we obtain $V_{2^n,p}=V_{4,p}$ for all $n\in\mathbb{N}_{\geq 2}$. This shows that the chain above stabilizes at $V_{4,p}$ and hence $V_{4,p}=\bar{k}$. \qed \\

\smallskip
\section{Proof of the Main Theorem}
\smallskip

In this section we will prove our main theorem. First, we require some additional tools. 
Let $\O$ be an order with conductor $\f$. Consider the multiplicative monoid \[\text{Reg}(\O)=\{ a\in \O \ | \  a+\f\in(\O/\f)^\times\}.\] It is easy to check that any element $a\in$ Reg$(\O)$ has the following property. If $x\in \O_K$, then $ax\in \O$ if and only if $x\in \O$.   \\ 

\noindent Note that there is a canonical embedding \[\O_K^\times/\O^\times\hookrightarrow (\O_K/\f)^\times/(\O/\f)^\times\] induced by restricting the natural surjection $\pi: \O_K\rightarrow\O_K/\f$ to $\O_K^\times$. We recall a fundamental result in the theory of orders in algebraic number fields.\\

\phantomsection
\label{4.1} \noindent \textbf{Lemma 4.1. (\cite[Theorem 2.11.5.2.]{HK20a})} \textit{We have} \[|\Pic(\O)|=s\cdot|\Pic(\O_K)|,\] \textit{where $s$ is the index of $\O_K^\times/\O^\times$ in $ (\O_K/\f)^\times/(\O/\f)^\times$.} \\

Recall that we call a prime ideal $\p$ of $\O$ \textit{inert} in $\O_K$ if $\p\O_K$ is again a prime ideal. \\

\phantomsection
\label{4.2} \noindent \textbf{Lemma 4.2.} \textit{We have $\O\cdot\O_K^\times=\O_K$ if and only if} $|\Pic(\O_K)|=|\Pic(\O)|$ \textit{and every prime ideal of $\O$ is inert in $\O_K$.} \\

\noindent \textit{Proof:} If $\O\cdot\O_K^\times=\O_K$, then we have an isomorphism \[\O_K^\times/\O^\times\cong (\O_K/\f)^\times/(\O/\f)^\times.\] Furthermore, for every prime ideal $\P$ of $\O_K$, we have $(\P\setminus\P^2)\cap\O\neq\emptyset$, which implies that no prime ideal of $\O$ is ramified in $\O_K$. It is left to check that the spec-map is bijective. Let $\P_1,\P_2$ be distinct prime ideals of $\O_K$, which lie over the prime ideal $\p$ and let $x\in\P_1\setminus\P_2$. If there is $\varepsilon\in\O_K^\times$ such that $y=\varepsilon x\in\O$, then $y\in\P_1\cap\O$, i.e. $y\in\p$, but this implies that $y\in\P_2$, which is a contradiction. \\

To prove the converse, suppose first that $\O$ is irreducible. Since every prime ideal of $\O$ is inert in $\O_K$, the spec-map is bijective and hence the conductor of $\O$ is of the form $\P^k$ for some positive integer $k$ and some prime ideal $\P$. Let now $a\in\P^r\setminus\P^{r+1}$ for some nonnegative integer $r$. Since $\P\cap\O$ is inert in $\O_K$, we find an element $b\in\O$ such that $b\in\P^r\setminus\P^{r+1}$. It is a well-known fact that one can find an element $u\in\O_K\setminus\P$ such that \[(u+\P^k)(a+\P^k)=b+\P^k \] and by Lemma \hyperref[4.2]{4.2.}, we can assume that $u\in \O_K^\times$. Hence we obtain that $a\O_K^\times\cap\O\neq\emptyset$. If $\O$ is not irreducible, then we obtain the same result by decomposing $\O/\f$ as a product of the residue class rings of the irreducible components of $\O$ with respect to their conductors and by applying the Chinese Remainder Theorem. \qed \\  

\phantomsection
\label{4.3} \noindent \textbf{Lemma 4.3.} \textit{Suppose that conditions (i) and (ii) of Theorem \hyperref[1.1]{1.1.} hold. Then $\O$ is half-factorial if and only if $\mathcal{A}(\O)\subseteq \mathcal{A}(\O_K)$.} \\

\noindent \textit{Proof:} If every atom of $\O$ is an atom of $\O_K$, then every factorization of a given element $x\in\O$ in $\O$ is also a factorization of $x$ in $\O_K$. Hence $\O$ is half-factorial. \\

Conversely, let $\O$ be half-factorial and assume that there is an element $x\in \O$, which is an atom of $\O$ but not an atom of $\O_K$. Then we can write \[x=\varepsilon a_1\cdots a_s,\] where $s\geq 2$, $\varepsilon\in \O_K^\times$ and every $a_i$ is an atom of $\O_K$. Since we have $\O\cdot \O_K^\times=\O_K$, we can assume that $a_i\in\O$ for every $i$. Then clearly $\varepsilon\not\in\O^\times$ and because $(\O_K^\times:\O^\times)$ is finite, we find a positive integer $n$, such that $\varepsilon^n\in\O^\times$. But then \[x^n=(\varepsilon^na_1) a_1^{n-1} a_2^n\cdots a_s^n\] and we see that $x^n$ has factorizations of length $n$ and $ns>n$ in $\O$, which is a contradiction. \qed   \\

Suppose that $\O$ is irreducible with conductor $\f=\P^k$ and consider the canonical homomorphism of monoids \[\varphi:\O^\bullet\longrightarrow(\O_\p^\bullet)_\text{red}\] \[x\mapsto x\O_\p^\times.\] Note that this map is surjective. Let $\varepsilon\in\O_K^\times$ and $x\in\O^\bullet$ such that $\varepsilon x\in\O$. Then, using the description of $(\O_\p^\bullet)_\text{red}$, which we gave in Section 3, we can write $\varphi(\varepsilon x)=u\varphi(x)$, where $u\in(\O_K/\f)^\times/(\O/\f)^\times$ corresponds to $\varepsilon\O^\times$ via the embedding \[\O_K^\times/\O^\times\hookrightarrow (\O_K/\f)^\times/(\O/\f)^\times.\] Suppose now that $\O\cdot\O_K^\times=\O_K$. Then for any two elements ${\bar{p}}^ru,{\bar{p}}^rv\in(\O_\p^\bullet)_\text{red}$ of the same value and $x\in\O$ with $\varphi(x)={\bar{p}}^ru$, we find by Lemma \hyperref[4.1]{4.1.} and Lemma \hyperref[4.2]{4.2.} a unit $\varepsilon\in\O_K^\times$ such that $\varphi(\varepsilon x)={\bar{p}}^rv$. \\

\noindent We will now prove one implication of our main theorem for irreducible orders.\\

\phantomsection
\label{4.4} \noindent \textbf{Proposition 4.4.} \textit{Let $\O$ be an irreducible order with conductor $\f=\P^k$ and suppose that conditions (i)-(iii) of Theorem \hyperref[1.1]{1.1.} hold. Then $\O$ is half-factorial.} \\

\noindent \textit{Proof:} By Lemma \hyperref[4.3]{4.3.}, we have to show that every atom of $\O$ is an atom of $\O_K$. Suppose to the contrary that there is $x\in\O$, which is an atom of $\O$ but not of $\O_K$. Let \[x\O_K=\P^\ell\mathfrak{a},\] where $\ell\in\mathbb{N}_0$ and $\mathfrak{a}$ is an ideal of $\O_K$, relatively prime to $\P$. If $\mathfrak{a}$ is contained in a proper principal ideal, say $\mathfrak{a}\subseteq\mathfrak{b}$, then we find elements $a\in \O_K$ and $b\in\O_K\setminus\P$ such that $b\O_K=\mathfrak{b}$ and $ab=x$. 
Since $\O\cdot\O_K^\times=\O_K$, we can assume that $b\in\O$. But then $b\in \text{Reg}(\O)$ and thus $a\in\O$. Since $x$ is an atom of $\O$, we obtain $a\in\O^\times$. This implies that $\ell=0$ and that $\mathfrak{a}$ is a principal ideal, which is not a product of two proper principal ideals. But then $x$ is an atom of $\O_K$, which is a contradiction to our assumption. Thus $\mathfrak{a}$ cannot be contained in a proper principal ideal and so we have to check three cases.  \\

\noindent \textbf{Case 1:} $\P$ is principal and $x\O_K=\P^\ell$ with $\ell\geq 2$. \\

Then $x$ has a unique factorization $\varepsilon\pi^\ell$ in $\O_K$, where $\pi$ is a prime element in $\O_K$. Let $\p=\P\cap\O$ and let $p\in\O_K$ be a prime element of $(\O_K)_\p$. Let $\varphi:\O^\bullet\rightarrow (\O_\p^\bullet)_{\text{red}}$ be the canonical homomorphism. Then $\mathsf{v}_{\bar{p}}(\varphi(x))=\ell$ and we can write $\varphi(x)=\bar{p}^\ell u$ for some $u\in (\O_K/\f)^\times/(\O/\f)^\times$, where $\bar{p}=p+\f$. By assumption, there is a factorization \[\bar{p}^\ell u=\bar{p}u_1\cdots \bar{p}u_\ell\] in $(\O_\p^\bullet)_{\text{red}}$. Since by Lemma \hyperref[4.1]{4.1.} and Lemma \hyperref[4.2]{4.2.}, we have an isomorphism \[\O_K^\times/\O^\times\cong (\O_K/\f)^\times/(\O/\f)^\times,\] 
we can choose $\varepsilon_1,\cdots,\varepsilon_\ell\in\O_K^\times$ with $\varphi(\varepsilon_i\pi)=\bar{p}u_i$ and hence $\varepsilon_i\pi\in\O$. But then \[\varepsilon^{-1}\varepsilon_1\cdots\varepsilon_\ell\in\O^\times\] and we see that \[\varepsilon_1\pi\cdots\varepsilon_\ell\pi\] is associate to $x$ in $\O$. This is a contradiction to $x$ being an atom of $\O$. \\

\noindent \textbf{Case 2:} $\P$ is not principal and $x\O_K=\P^{2m}$ with $m\geq 2$. \\

Then $x$ has a unique factorization $\varepsilon a^m$ in $\O_K$, where $a$ is an atom of $\O_K$. We have $\varphi(x)=\bar{p}^{2m}u$ and by assumption, there is a factorization \[\bar{p}^{2m}u=\bar{p}u_1\cdots \bar{p}u_r \bar{p}^2u_{r+1}\cdots \bar{p}^2u_s.\] Writing \[\bar{p}^{2m}u=\bar{p}^2(u_1u_2)\cdots \bar{p}^2(u_{r-1}u_r) \bar{p}^2u_{r+1}\cdots \bar{p}^2u_s,\] we can again find $\varepsilon_1,\cdots,\varepsilon_m\in\O_K^\times$ such that the $\varepsilon_i a$-s correspond to these factors via $\varphi$. But then $\varepsilon_ia\in\O$ and \[\varepsilon_1a\cdots\varepsilon_ma\] is associate to $x$ in $\O$.   \\

\noindent \textbf{Case 3:} $\P$ is not principal and $x\O_K=\P^{2m+1}\mathfrak{Q}$, where $m\in\mathbb{N}_0$ and $\mathfrak{Q}$ is a non-principal prime ideal of $\O_K$. \\

Then $x$ has a unique factorization $\varepsilon a^m b$ in $\O_K$, where $a$ and $b$ are atoms of $\O_K$ such that $a\O_K=\P^2$ and $b\O_K=\P\mathfrak{Q}$. This case is then analogous to Case 2. \qed \\

\noindent \textit{Proof of Theorem \hyperref[1.1]{1.1.}:} Let $\O$ be half-factorial. We have to show that conditions (i)-(iii) of Theorem \hyperref[1.1]{1.1.} hold. \\

By \cite[Example 3.7.3]{Ge-HK06a} or \cite[Corollary 4.6.]{Ph12b}, we see that $|\Pic(\O)|\leq 2$. By Lemma \hyperref[4.1]{4.1.}, we then have $|\Pic(\O_K)|\leq 2$ and $\O_K$ is half-factorial. \\

Since the elasticity $\rho(\O)$ is finite, we obtain by  \cite[Corollary 3.7.2]{Ge-HK06a} that the spec-map is bijective. \\

Suppose that $\O_K$ is a PID and that $\O\cdot\O_K^\times\neq\O_K$. Then there is a prime element $\pi$ with $\pi\O_K^\times\cap\O=\emptyset$. Take an element $f\in\f$. Then $\pi f\in\O$ and by Lemma \hyperref[2.1]{2.1.}, we find $n\geq 2$ such that $\pi^n\in\O$. But then both $(\pi f)^n$ and $f^n$ have factorizations in $\O$ for which the lengths are multiples of $n$. This implies that $\pi^n$ has a factorization for which the length is a multiple of $n$. This is impossible, since there is no $\varepsilon\in\O_K^\times$ with $\varepsilon\pi\in\O$. So we have shown that $\O\cdot\O_K^\times=\O_K$. \\

Let $\P$ be a prime ideal of $\O_K$ with $\P\supseteq\f$ and let $\pi$ be a prime element such that $\pi\O_K=\P$. We can assume that $\pi\in\O$. We claim that no element in $\P^2\cap\O$ is an atom of $\O$. Assume the contrary, i.e. there is an element $\pi^ka\in \O$ with $k\geq 2$, which is an atom of $\O$. Clearly, $a\not\in\O$. Then by Lemma \hyperref[2.1]{2.1.}, we find $n\geq 2$ such that $a^n\in\O$. Then $(\pi^k a)^n$ has a factorization of length $n$ in $\O$, but $\pi^{kn}$ has a factorization of length $kn>n$, which is a contradiction. Let $\p=\P\cap \O$ and let $p$ be a prime element of $(\O_K)_\p$. Note that any factorization of an element $x\in\P^2\cap\O$ induces a factorization of the corresponding element $\varphi(x) \in (\O_\p^\bullet)_{\text{red}}$. Since $\varphi$ is surjective, we obtain $\mathsf{v}_p(\mathcal{A}(\O_\p))=\{1\}$ and thus $\O_\p$ is half-factorial.  \\

\noindent This proves conditions (i)-(iii) in the case that $\O_K$ is a PID. \\

If $\O_K$ is not a PID, then $|\Pic(\O_K)|=|\Pic(\O)|=2$ and by Lemma \hyperref[4.1]{4.1.}, we have an isomorphism \[\O_K^\times/\O^\times\cong (\O_K/\f)^\times/(\O/\f)^\times.\] Let $\pi$ be a prime element in $\O_K$. Then by the argument above, we see that $\pi\O_K^\times\cap\O\neq\emptyset$. The same holds for an atom $a$ of $\O_K$, for which $a\O_K=\P^2$ for a non-principal prime ideal $\P$, since every power of $a$ has a unique factorization in $\O_K$. \\

Let $\P$ be principal with $\P\supseteq \f$. By repeating the argument above, we see that no element in $\P^2\cap\O$ can be an atom of $\O$, which implies condition (iii) for principal ideals.  \\

Let now $\P\supseteq\f$ be non-principal and suppose that there is an element $x\in\P^3\cap\O$, which is an atom of $\O$. Then we can write $x=ab$, where $a\O_K=\P^2$ and $b\in\P$. As mentioned above, we can assume that $a\in\O$. Then, clearly, $b\not\in\O$. Again, by Lemma \hyperref[2.1]{2.1.}, we find $n\geq 2$ such that $b^n\in\O$. But then $x^n=(ab)^n$ and $a^n$ both have factorizations of length $n$, which implies that $b^n$ has a factorization of length $0$, i.e. $b^n$ is a unit. Since $b\in\P$, this is clearly impossible. This implies that $\mathsf{v}_p(\mathcal{A}(\O_\p))\subseteq\{1,2\}$. Clearly, $\mathsf{v}_p(\mathcal{A}(\O_\p))=\{2\}$ is impossible, since $\O_\p$ has finite exponent. \\

So every localization contains an element of value $1$, which shows that every prime ideal of $\O$ is inert in $\O_K$. By Lemma \hyperref[4.2]{4.2.}, we obtain $\O\cdot\O_K^\times=\O_K$ and Propositions \hyperref[3.2]{3.2.} and \hyperref[3.3]{3.3.} yield the bounds on the exponents of the conductor. Thus we have shown all of the conditions (i)-(iii). \\

For the converse, suppose that the conditions (i)-(iii) of Theorem \hyperref[1.1]{1.1.} hold. We have to show that $\O$ is half-factorial. By Lemma \hyperref[4.2]{4.2.}, we see that every prime ideal of $\O$ is inert in $\O_K$, which implies that the spec-map is bijective and so we can write \[\O=\bigcap_{i=1}^{s}\O_i,\] where the $\O_i$'s are the irreducible components of $\O$ with conductor $\P_i^{k_i}$ respectively. Clearly, conditions (i) and (ii) hold for every $\O_i$ and by the isomorphism given in Lemma \hyperref[3.1]{3.1.}, condition (iii) holds as well. Hence by Proposition \hyperref[4.4]{4.4.}, every $\O_i$ is half-factorial. \\

By Lemma \hyperref[4.3]{4.3.}, it is enough to show that every atom of $\O$ is an atom of $\O_K$. Assume to the contrary that there is $x\in\O$, which is an atom of $\O$ but not of $\O_K$. Clearly, $x$ can not be divided by an element in $\text{Reg}(\O)\setminus\O^\times$, which implies that $x\O_K$ is of the form \[\P_1^{r_1}\cdots\P_s^{r_s}\mathfrak{Q}^r,\] where $r\in\{0,1\}$, $r_i\in \mathbb{N}_0$ and $\mathfrak{Q}$ is a non-principal prime ideal, relatively prime to $\f$. \\

Take $i\in [1,s]$ such that $r_i\neq 0$ and supose that $\P_i$ is principal. Let $\pi$ be a prime element such that $\pi\O_K=\P$ and assume that $\pi\in\O$. We can write $x=\varepsilon\pi_i^{r_i}y$ for some $\varepsilon\in\O_K^\times$ and some $y\in\O\setminus\P_i$.
Since we have isomorphisms \[\O_K^\times/\O^\times\cong (\O_K/\f)^\times/(\O/\f)^\times\cong \bigtimes_{\ell=1}^{s}(\O_K/\P_{\ell}^{k_\ell})^\times /(\O_\ell/\P_\ell^{k_\ell})^\times,\] we can write $\varepsilon\O^\times=\begin{pmatrix}
    u_1,\cdots,u_s
\end{pmatrix}$. There is $\eta\in\O_K^\times$ with \[\eta\O^\times=\begin{pmatrix}
    1,\cdots,u_i,\cdots,1
\end{pmatrix}.\] Then $\eta\pi_i^{r_i}\in\O$ and $\eta^{-1}\varepsilon y\in\O$, and since $x$ is an atom of $\O$, we see that $y\in\O^\times$ and thus $x\O_K=\P^{r_i}$. \\

The same type of argument can be applied in the following cases. 

\begin{enumerate}
    \item[$\bullet$] $\P_i$ is non-principal with $r_i\neq 0$ and $r_i$ even. Simply replace $\pi_i$ with an atom $a_i\in\O$ such that $a_i\O_K=\P_i^2$.
    \item[$\bullet$] $\P_i$ is non-principal with $r=1$, $r_i\neq 0$ and $r_i$ odd. Then we have atoms $a_i\in \O$ and $b\in \O$ such that $a\O_K=\P_i^2$ and $b\O_K=\P_i\mathfrak{Q}$ and we replace $\pi_i^{r_i}$ with $a_i^mb$, where $r_i=2m+1$.
    \item[$\bullet$] $\P_i$ and $\P_j$ are non-principal with $r_i,r_j\neq 0$ and odd. Then we find an element $z\in\O$ with $z\O_K=\P_i^{r_i}\P_j^{r_j}$ such that $x=\varepsilon zy$ for some $\varepsilon\in\O_K^\times$ and $y\in\O\setminus(\P_i\cup\P_j)$. Then, if $\varepsilon\O^\times=\begin{pmatrix}
        u_1,\cdots,u_s
    \end{pmatrix}$, we find an $\eta\in\O_K^\times$ with \[\eta\O^\times=\begin{pmatrix}
        1,\cdots,u_i,\cdots,1,\cdots,u_j,\cdots,1
    \end{pmatrix}\] and we see again that $y\in\O^\times$.
\end{enumerate}

Hence we are left with four cases. \\

  \noindent \textbf{Case 1:} $x\O_K=\P_i^{r_i}$, where $i\in[1,s]$, $\P_i$ is principal, and $r_i\geq 2$. \\
    
    Thus $x=\varepsilon \pi_i^{r_i}$ for some $\varepsilon\in\O_K^\times$. Since $\O_i$ is half-factorial, there is a factorization \[x=\varepsilon_1\pi_i\cdots\varepsilon_{r_i}\pi_i\] in $\O_i$. For $j\in[1,r_i]$, let $\varepsilon_j\O^\times=\begin{pmatrix}
        u_{1,j},\cdots,u_{s,j}
    \end{pmatrix}$. We can find $\eta_j\in\O_K^\times$ with $\eta_j\O^\times=\begin{pmatrix}
        1,\cdots,u_{i,j},\cdots,1
    \end{pmatrix}$. Then $\eta_j\pi_i\in\O$ for all $j\in[1,r_i]$ and \[\eta_1\pi_i\cdots\eta_{r_i}\pi_i\] is associate to $x$ in $\O$. Hence $x$ is not an atom of $\O$, which is a contradiction to our assumption. \\
    
   \noindent \textbf{Case 2:} $x\O_K=\P_i^{r_i}$, where $i\in[1,s]$, $\P_i$ is non-principal, and $r_i\geq 4$ is even. \\

    Repeating the argument given in the first case by replacing $\pi_i$ with an atom $a_i$ with the property $a_i\O_K=\P_i^2$ yields that $x$ is not an atom of $\O$, which is again a contradiction to our assumption. \\

 \noindent \textbf{Case 3:} $x\O_K=\P_i^{2r_i+1}\mathfrak{Q}$, where $i\in[1,s]$, $\P_i$ and $\mathfrak{Q}$ are non-principal, $\mathfrak{Q}$ is relatively prime to $\f$ and $r_i\in\mathbb{N}$. \\

  Let $a,b\in\O$ such that  $a\O_K={\P_i}^2$ and $b\O_K=\P_i\mathfrak{Q}$. By replacing $\pi_i^{r_i}$ with $a_i^{r_i}b_i$ and proceeding as in the first case, it follows that $x$ is not an atom of $\O$, another contradiction to our assumption. \\

\noindent \noindent \textbf{Case 4:} $x\O_K=\P_i^{r_i}\P_j^{r_j}$, where $i,j\in[1,s]$ with $i\neq j$; $\P_i,\P_j$ are non-principal, and $r_i=2m_i+1,r_j=2m_j+1$ are both odd with $m_i,m_j\in\mathbb{N}_0$ and $m_i+m_j\geq 1$. \\

Let $a_i,a_j,b\in\O$ such that $a_i\O_K=\P_i^2$, $a_j\O_K=\P_j^2$ and $b\O_K=\P_i\P_j$ and write \[x=\varepsilon a_i^{m_i}a_j^{m_j}b.\] Note that $a_i\in\text{Reg}(\O_j)$ and $a_j\in\text{Reg}(\O_i)$ and so $\varepsilon a_i^{m_i}b\in\O_i$ and $\varepsilon a_j^{m_j}b\in\O_j$. Since both $\O_i$ and $\O_j$ are half-factorial, we obtain the following factorizations of $x$. \[x=a_j^{m_j}\cdot \eta b\cdot\prod_{l=1}^{m_i}\eta_l a_i \text{ in } \O_i\] and \[x=a_i^{m_i}\cdot \xi b\cdot\prod_{\ell=1}^{m_j}\xi_\ell a_j \text{ in } \O_j\] for units $\eta,\eta_l,\xi,\xi_\ell\in \O_K^\times$, where $l\in[1,m_i]$ and $\ell\in[1,m_j]$. Let \[\eta\O^\times=\begin{pmatrix}
        \tilde{\eta}_1,\cdots,\tilde{\eta_s}
    \end{pmatrix},\] \[\eta_l\O^\times=\begin{pmatrix}
        \tilde{\eta}_{l,1},\cdots,\tilde{\eta}_{l,s}
    \end{pmatrix},\]\[\xi\O^\times=\begin{pmatrix}
        \tilde{\xi}_1,\cdots,\tilde{\xi}_s
    \end{pmatrix},\]\[\xi_\ell\O^\times=\begin{pmatrix}
        \tilde{\xi}_{\ell,1},\cdots,\tilde{\xi}_{\ell,s}
    \end{pmatrix}.\] Then we can find units $\mu_l,\nu_\ell,\zeta\in\O_K^\times$ with \[\mu_l\O^\times=\begin{pmatrix}
        1,\cdots,\tilde{\eta}_{l,i},\cdots,1
    \end{pmatrix},\] \[\nu_\ell\O^\times=\begin{pmatrix}
        1,\cdots,\tilde{\xi}_{\ell,j},\cdots,1
    \end{pmatrix},\] and \[\zeta\O^\times=\begin{pmatrix}
        1,\cdots,\tilde{\eta}_i,\cdots,1,\cdots,\tilde{\xi}_j,\cdots,1
    \end{pmatrix}.\] Thus $\mu_la_i, \nu_\ell a_j$ and $\zeta b$ are contained in $\O$ for all $l\in[1,m_i], \ell\in[1,m_j]$ and \[\zeta b\cdot\prod_{l=1}^{m_i}\mu_la_i\cdot\prod_{\ell=1}^{m_j}\nu_\ell a_j\] is associate to $x$ in $\O$, which yields the final contradiction. \qed \\
    
We conclude with two corollaries and some examples of half-factorial orders. The first corollary is the characterization of half-factorial seminormal orders, proved by Geroldinger, Kainrath and Reinhart in a more general setting in \cite[Theorem 6.2]{Ge-Ka-Re15a}. \\

A Noetherian domain $R$ is called \textit{seminormal} if there is no $x\in\bar{R}\setminus R$, such that $x^n\in R$ for almost all $n\in\mathbb{N}$, where $\bar{R}$ is the integral closure of $R$. \\

\phantomsection
\label{4.5} \noindent \textbf{Proposition 4.5.} \textit{Let $K$ be an algebraic number field and let $\O$ be an order in $K$ with conductor $\f$. Then $\O$ is seminormal if and only if $\f$ is squarefree.} \\

\noindent \textit{Proof:} Suppose that $\f$ is squarefree and let $\O_i$ be an irreducible component of $\O$ with conductor $\f_i$. If there is $x\in\O_K$ such that $x^n\in\O$ for almost all $n\in\mathbb{N}$, then $x+\f_i$ is clearly either nilpotent or a unit in $\O_K/\f_i$. Since $\O_K/\f_i$ is a product of fields, we have in both cases $x\in\O_i$ and hence $x\in\O$.  \\

If $\f$ is not squarefree, then there is a prime ideal $\P$ of $\O_K$, such that $\P^2$ divides $\f$. Any element $x \in \f\P^{-1}$ satisfies $x^n\in\O$ for almost all $n$, however this ideal is not contained in $\O$ and thus $\O$ is not seminormal. \qed\\

\phantomsection
\label{4.6} \noindent \textbf{Corollary 4.6.} \textit{Let $K$ be an algebraic number field and let $\O$ be a seminormal order in $K$. Then $\O$ is half-factorial if and only if the following conditions are satisfied } \\

(i) \textit{$\O_K$ is half-factorial.} \\

(ii) \textit{The spec-map is bijective.} \\

(iii) $|\Pic(\O_K)|=|\Pic(\O)|$. \\

\noindent \textit{Proof}: If $\O$ is half-factorial, all the conditions follow of Theorem \hyperref[1.1]{1.1.} and Lemma \hyperref[4.2]{4.2.}. \\

Conversely, it follows from Proposition \hyperref[4.5]{4.5.} that every localization satisfies $\mathsf{v}_p(\mathcal{A}(\O_\p))=\{1\}$. Since the spec-map is bijective, Lemma \hyperref[4.2]{4.2.} then implies that $\O\cdot\O_K^\times=\O_K$ and by Theorem \hyperref[1.1]{1.1.}, $\O$ is half-factorial. \qed \\
 
Our second corollary concerns the characterization of half-factorial orders in quadratic number fields, which was proved by Halter-Koch, see \cite[§3, Theorem 6]{HK83a} or \cite[Theorem 3.7.15.]{Ge-HK06a}. \\

First, we recall some properties of quadratic orders. For a detailed background, we refer to \cite{HK13a}. Let $K$ be a quadratic number field. Every conductor ideal is of the form $f\O_K$ for some $f\in\mathbb{N}$ and we will just say that an order $\O$ has conductor $f$. Furthermore, the only order with conductor $f$ is the minimal order $\mathbb{Z}+f\O_K$.   \\

By Dirichlet's Unit Theorem, we have $\O_K^\times\cong \{\pm 1\}\times \mathbb{Z}$ if $K$ is real quadratic and $\O_K^\times$ is finite of order at most $6$ if $K$ is imaginary quadratic. Since we have $\{\pm 1\}\subseteq \O^\times$ for any order $\O$, we obtain that $\O_K^\times/\O^\times$ is either isomorphic to $\mathbb{Z}$ or finite of order at most 3 and hence cyclic in both cases.\\

\phantomsection
\label{4.7} \noindent \textbf{Corollary 4.7.} \textit{Let $\O$ be an order in a quadratic number field $K$ with conductor $f\in\mathbb{N}_{\geq 2}$. Then $\O$ is half-factorial if and only if the following conditions are satisfied} \\

(i) \textit{$\O_K$ is half-factorial.} \\

(ii) \textit{$\O\cdot\O_K^\times=\O_K$.} \\

(iii) \textit{$f$ is either a prime or twice an odd prime.} \\

\noindent \textit{Proof:} Let $\O$ be half-factorial. By Theorem \hyperref[1.1]{1.1.}, we only have to prove that $f$ is a prime or twice an odd prime. Assume first that $\O$ is irreducible. Since the spec-map is bijective, the conductor of $\O$ is of the form $f\O_K=\P^k$, i.e. $f=p^r$ for a prime integer $p$ that does not split in $K$. \\

If $p$ is ramified in $K$, i.e. $p\O_K=\P^2$, we have $\O=\mathbb{Z}+\P^{2r}$. Note that no integer can be contained in $\P\setminus\P^2$ and so $\P\cap\O$ is not inert in $\O_K$. Hence by Lemma \hyperref[4.2]{4.2.}, the condition $\O\cdot\O_K^\times=\O_K$ is violated. \\

If $p$ is inert in $K$, the prime ideal $p\O_K=\P$ is principal and by Theorem \hyperref[1.1]{1.1.}, we have $r\leq 2$. If $r=2$ and $\p=\P\cap\O$, we see by Proposition \hyperref[3.2]{3.2.} (c) and the fact that $p$ has inertia degree two, that $\O_\p$ is not half-factorial and thus $r=1$. \\

Hence, if $\O$ is not irreducible, then $f$ is a squarefree product of inert primes. Suppose that there are two odd primes $p,q$, which divide $f$. Since \[|(\O_K/p\O_K)^\times/((\mathbb{Z}+p\O_K)/p\O_K)^\times|=p+1\] and \[|(\O_K/q\O_K)^\times/((\mathbb{Z}+q\O_K)/q\O_K)^\times|=q+1,\] we see that $(\O_K/f\O_K)^\times/(\O/f\O_K)^\times$ is an abelian group of rank at least two. Since $\O_K^\times/\O^\times$ is cyclic, we obtain a contradiction. \\

Conversely, let $\O_K$ be half-factorial, let $\O\cdot\O_K^\times=\O_K$ and let $f$ be a prime or twice an odd prime. Repeating the argument from above, we easily see that no ramified prime can divide $f$. Similarly, if $p$ divides $f$, it cannot split in $K$, since the spec-map is bijective by Lemma \hyperref[4.2]{4.2.}. Hence every prime that divides $f$ is inert, and since $f$ is square-free, the ideal $f\O_K$ is squarefree as well. This implies that all localizations $\O_\p$ have exponent one and therefore, they are half-factorial. Then by Theorem \hyperref[1.1]{1.1.}, $\O$ is half-factorial. \qed  \\

In \cite{Ka05b}, Kainrath indirectly conjectured that every order between a half-factorial order and $\O_K$ is itself half-factorial. Although our main theorem suggests that a counterexample could be constructed, this conjecture remains open.  \\

To conclude, we will present some examples of half-factorial orders and show that none of the conditions in Theorem 1.1. is vacuous.\\

\text{ }\\

\noindent \textbf{Examples:} \\

\noindent (1) In \cite{Br-Ge-Re20}, the authors investigate the arithmetic of orders in quadratic number fields. In particular, \cite[Theorem 1.1.]{Br-Ge-Re20} gives a characterization of the half-factoriality of the monoid of invertible ideals $\mathcal{I}^\ast(\O)$ of a quadratic order $\O$. Similar results about stable domains can be found in \cite{Ba-Ge-Re21c}. If $|\Pic(\O)|=1$, then $\mathcal{I}^\ast(\O)$ being half-factorial is equivalent to $\O$ being half-factorial. In general, for an order $\O$, we have \[\mathcal{I}^\ast(\O)\cong \coprod_{\p\in\spec(\O)}(\O_\p^\bullet)_{\text{red}}\] and hence $\mathcal{I}^\ast(\O)$ is half-factorial if and only if every localization is half-factorial. In (3), we will see that if $|\Pic(\O)|=2$, then the half-factoriality of $\O$ does not imply the half-factoriality of $\mathcal{I}^\ast(\O)$.\\

\noindent (2) In \cite{Po-24}, the author shows that there are infinitely many half-factorial real quadratic orders. Moreover, under the assumption of the generalized Riemann hypothesis, $\mathbb{Q}(\sqrt{2})$ contains infinitely many half-factorial orders. These results are in stark contrast to the case of imaginary quadratic orders, since in \cite[Theorem 7]{Ch-Co00}, it is shown that $\mathbb{Z}[\sqrt{-3}]$ is the only half-factorial imaginary quadratic (non-principal) order.  \\

\noindent (3) We will now give an example of a half-factorial order $\O$, for which there is a prime ideal $\p\subseteq \O$ such that $\mathsf{v}_p(\mathcal{A}(\O_\p))=\{1,2\}$. In other words, we prove the existence of a half-factorial order, of which not every localization is half-factorial. In particular, $\mathcal{I}^\ast(\O)$ is not half-factorial. \\

Let $f=X^3-8X-19\in\mathbb{Z}[X]$, let $\omega$ be a root of $f$ and set $K=\mathbb{Q}(\omega)$. We gather the following facts from \cite[Number field 3.1.7699.1]{LM-24} and from computations with SageMath.\\
\begin{itemize}
    \item We have $\O_K=\mathbb{Z}[\omega]$ and $|\Pic(\O_K)|=2$.
    \item $\O_K^\times\cong \mathbb{Z}\times\{\pm 1\}$ with $\varepsilon=15\omega^2-32\omega-82$ being a fundamental unit.
    \item The principal ideal $2\O_K$ can be factored as \[2\O_K=(2,\omega^2-\omega-5)\cdot(2,\omega+1).\]

\end{itemize}

\noindent Let $\P=(2,\omega^2-\omega-5)$ and let $\f=\P^2=(\omega^2-5\omega+5)$. We set $\O:=\mathbb{Z}+\f$. Since $\P$ has inertia degree $2$, we see that $|\O_K/\f|=16$. Identifying $(\O_K/\f,+)$ with $(\mathbb{Z}/4\mathbb{Z})^2$ and writing \[\O_K/\f=\{a\omega+b:a,b\in\mathbb{Z}/4\mathbb{Z}\},\] we obtain that \[(\O_K/\f)^\times=(\O_K/\f)\setminus\{\overline{0},\overline{2},\overline{2}\omega,\overline{2}\omega+\overline{2}\}.\] Moreover, we have \[\O/\f=\{\overline{0},\overline{1},\overline{2},\overline{3}\}\] and \[(\O/\f)^\times=\{\overline{1},\overline{3}\}.\] 
A quick computation yields that $(\O_K^\times:\O^\times)=6$ and by Lemma \hyperref[4.1]{4.1.,} we obtain that $|\Pic(\O_K)|=|\Pic(\O)|$. Since $\p=\P\cap\O$ is clearly inert in $\O_K$, we also have $\O\cdot\O_K^\times=\O_K$ by Lemma \hyperref[4.2]{4.2.} It is easy to see that \[\{\overline{1},\overline{3},\overline{2}\omega+\overline{1},\overline{2}\omega+\overline{3}\}\] is the set of all elements $\varepsilon\in (\O_K/\f)^\times$ such that $\overline{2}\cdot\varepsilon\in \O/\f$. Using the same notation as in the remark preceding Lemma \hyperref[3.1]{3.1.}, we obtain that \[\overline{U}_{1,\bar{p}}=\{\overline{1}(\O/\f)^\times,(\overline{2}\omega+\overline{1})(\O/\f)^\times\},\] where $p=\overline{2}$. 
However, $\overline{U}_{1,\bar{p}}$ is a group of order 2 under multiplication and hence \[\overline{U}_{1,\bar{p}}\cdot\overline{U}_{1,\bar{p}}\neq \overline{U}_{2,\bar{p}}=(\O_K/\f)^\times/(\O/\f)^\times.\] Consequently, $(\O_\p^\bullet)_{\text{red}}$ contains an atom of value 2 but no atoms of value 3 or higher. Then, by Theorem \hyperref[1.1]{1.1.}, $\O$ is a half-factorial order with $\mathsf{v}_q(\mathcal{A}(\O_\p))=\{1,2\}$, where $q$ is some prime element of $(\O_K)_\p$. \\

\noindent \textbf{Acknowledgements} I would like to thank Alfred Geroldinger and the anonymous referee for their helpful comments and suggestions.

\providecommand{\bysame}{\leavevmode\hbox to3em{\hrulefill}\thinspace}
\providecommand{\MR}{\relax\ifhmode\unskip\space\fi MR }
\providecommand{\MRhref}[2]{%
  \href{http://www.ams.org/mathscinet-getitem?mr=#1}{#2}
}
\providecommand{\href}[2]{#2}

\end{document}